\documentclass[oneside]{amsart}

\pdfoutput=1
\usepackage[paper, final, normal, externalize=false, journal=ams]{style}

	\author[]{Sebastian Halbig} 
	\address[]{
		Technische Universität Dresden,
		Institut f\"ur Geometrie,
		\\
		Zellescher Weg 12--14, 01062 Dresden}
	\email[]{sebastian.halbig@tu-dresden.de}
	\date{\today}
	\title{Generalised Taft algebras and many pairs in involution}
\begin{document}	
	\numberwithin{theorem}{section}
	\begin{abstract}
		A class of finite-dimensional Hopf algebras which generalise the notion of Taft algebras is studied. 
		We give necessary and sufficient conditions for these Hopf algebras to 
		omit a pair in involution.
		That is, to not have a group-like and a character implementing the square of the antipode.
		As a consequence, we prove the existence of an infinite set of examples of finite-dimensional Hopf algebras without such pairs.
		Implications for the theory of anti-\YetterDrinfeld modules as well as biduality of representations of Hopf algebras are discussed.
	\end{abstract}

	\maketitle
	
	\paragraph{Acknowledgements}	
		The author would like to thank P. Hajac for his kind invitation to IMPAN. He would also like to thank I. Heckenberger and 
		U. Kr\"ahmer for many stimulating discussions.
		
	\section{Introduction}
	\subsection{Main result}
	A pair in involution for a Hopf algebra $H$ over a field $\k$  is a pair $(l, \beta)$ of a group-like  $l\in H$ and a character $\beta \from H \to  \k$ such that the square of the antipode is given by the conjugate action of $l$ and $\beta$. 
	In a vague sense, it can be  imagined to be similar to a `square-root' of Radford's $S^4$-formula. 
	Often one additionally requires the pair to be modular, i.e to satisfy $\beta(l) = 1$.
	Pairs in involution appear in many different contexts within Hopf algebra theory,
	reaching from Hopf-cyclic cohomology \cite{Connes1999, Hajac2004, Hajac2010} to knot invariants \cite{Kauffman1993}.
	Kauffman and Radford showed in the aforementioned article that the square of the antipode of certain Taft algebras is not implemented by square-roots of the distinguished group-likes, see \cite[Proposition 7]{Kauffman1993}.
	Nonetheless, these Hopf algebras admit modular pairs in involution.
	In a previous work by Kr\"ahmer and the author, see \cite{Halbig2019}, examples of Hopf algebras were given, whose pairs in involution do not satisfy the modularity condition.
	The paper at hand builds upon this result and proves the existence of  finite-dimensional Hopf algebras without such pairs.		
	Hereto, we introduce generalised Taft algebras, a class of Hopf algebras containing the examples of \cite{Halbig2019}.
	As algebras, these are generated by a group-like $g$ and two twisted primitives $x$ and $y$ such that $g$ spans a cyclic group and $x$ and $y$ are nilpotent. Moreover, the generators are required to commute up to some 
	roots of unity. 
	The details of this presentation are given in Theorem \ref{thm: GTAGenerators}.
	For a generalised Taft algebras pairs in involution correspond to solutions of systems of Diophantine equations, see Theorem \ref{thm: ModPairAsEquation}.
	The main result, Theorem \ref{thm: Classification}, gives necessary and sufficient conditions for 
	the non-existence of such solutions.
	In Lemma \ref{lem: Examples} we apply the above result to show that
	there are finite-dimensional Hopf algebras without these pairs.
	Additionally, Lemma \ref{lem: IsoAlgNoPair} shows that Hopf algebras which are isomorphic only as algebras need not share the property of having a pair in involution.	
	\subsection{Pivotal categories and anti-\YetterDrinfeld modules} 

	A pivotal category is a monoidal category with a notion of duality and a natural  isomorphism between any object and its bidual which is compatible with the monoidal structure.
	An example are the finite-dimensional modules over pivotal Hopf algebras. 
	As discussed in Section \ref{sec: ADDandPiv}, a finite-dimensional Hopf algebra has a pair in involution if and only if its Drinfeld double is pivotal.
	Another way of describing this interplay is via anti-\YetterDrinfeld modules, which arose in the field of Hopf-cyclic cohomology. 
	Similar to \YetterDrinfeld modules, they are simultaneously modules and comodules  over a Hopf algebra with a compatibility condition between the action and coaction.
	For a finite-dimensional Hopf algebra one can construct an algebra called the anti-Drinfeld double, whose modules correspond to anti-\YetterDrinfeld modules.
	The existence of a pair in involution is equivalent to the Drinfeld double and the anti-Drinfeld double being isomorphic as algebras.
	As a consequence of our findings there are examples where such an isomorphism does not exist. In this case the Drinfeld double is not a pivotal Hopf algebra.
	
	\subsection{Outline}
	This article is organized as follows.
	Section \ref{sec: PII} serves as a summary of the theory of pairs in involution with a focus on finite-dimensional Hopf algebras.
	In Section \ref{sec: GTA} we introduce the main object of our study, generalised Taft algebras, and give an `easy-to-work-with' presentation.
	The classification of all of these Hopf algebras omitting a pair in involution is carried out in Section \ref{sec: HAWithoutPII}. 
	Afterwards we discuss examples of Hopf algberas with and without such pairs.
	We conclude the paper with Section \ref{sec: ADDandPiv}, where we apply our results to the context of representation theory.
		
	\section{Pairs in involution}\label{sec: PII}
	We work over an algebraically closed field $\k$ of characteristic zero; `$\dim$' and `$\otimes$' ought to be understood as dimension and tensor product over $\k$. 
	Standard notation for Hopf algebras, as in e.g.\ \cite{Montgomery1993,Radford2012}, is freely used. 
	Given a Hopf algebra $H$ we write $\Gr(H)$ for its group of group-likes, $\Pr(H)$ for its space of primitive elements and $H^\circ$ for its (finite) dual Hopf algebra. 
	The antipode of $H$ is denoted by   $S \from H \rightarrow H$, its counit by $\epsilon \from H \rightarrow \k$ and its coproduct by $\Delta \from H \rightarrow H\otimes H$. 
	For calculations involving the coproduct of $H$ or the coaction of some comodule $M$ over $H$ we rely on reduced Sweedler notation. 
	For example we  write $h_{(1)} \otimes h_{(2)} \defeq \Delta (h)$ for $h \in H$.
	An element $x\in H$ whose coproduct is $\Delta(x) = 1 \otimes x + x \otimes g$, for $g\in \Gr(H)$ a group-like, is called a twisted primitive.
	\newline
	
	Modular pairs in involution play the role of coefficients for Hopf-cyclic cohomology as introduced by Connes and Moscovici \cite{Connes1999}. 
	Later on, it was realised by Hajac \etal \cite{Hajac2004a} that this notion can be extended to that of (stable) anti-\YetterDrinfeld modules, which we discuss in Section \ref{sec: GTA}.
	
	\begin{definition}\label{def: PII}
		Let $H$ be a Hopf algebra over $\k$. A pair $(l, \beta )$ comprising group-like elements $l \in \Gr(H)$ and $\beta \in \Gr(H^\circ)$ is a \emph{pair in involution} if it satisfies the antipode condition 
		\begin{equation}\label{eq: AntipodeCondition}
			S^2(h) = \beta (h_{(3)}) \beta^{-1} (h_{(1)}) \;  l h_{(2)} l^{-1}
			,
			\qquad \text{ for all } h \in H. \tag{AC}
		\end{equation}
		If additionally the modularity condition
		\begin{equation}\label{eq: modularityCondition}
			\beta(l) = 1 \tag{MC}
		\end{equation}
		holds, it is called a \emph{modular} pair in involution.
	\end{definition} 
	
	A \emph{left integral} of a Hopf algebra $H$ is an element $\Lambda \in H$ such that 
	$h \Lambda  = \epsilon(h) \Lambda$ for all $h\in H$. 
	If $H$ is finite-dimensional, its left integrals form a one-dimensional subspace $\LInt(H)\subset H$. 
	There is a unique group-like  $\alpha\in \Gr(H^\circ)$ such that $\Lambda h= \alpha(h) \Lambda$ for all  $\Lambda \in \LInt(H)$ and $h \in H$. 
	It is called the \emph{distinguished group-like} of $H$.	
	Radford proved that the fourth power of the antipode is implemented by the distinguished group-likes of a Hopf algebra and its dual \cite{Radford1976}. 
	
	\begin{theorem}[Radford's $S^4$ formula]\label{thm: S4}
		Let $H$ be a finite-dimensional Hopf algebra and $g \in \Gr(H)$, $\alpha \in \Gr(H^\circ)$ the distinguished group-likes of $H^\circ$ and $H$. Then the fourth power of the antipode is given by
		\begin{equation} \label{eq: S4-formula}
			S^4(h)  
			= \alpha(h_{3})\alpha^{-1}(h_{(1)}) \; g^{-1}h_{(2)}g, 
			\quad \text {for all } h \in H. 
		\end{equation} 
	\end{theorem}

	In their paper on the classification of ribbon elements of Drinfeld doubles Kauffman and Radford studied `square roots' of the distinguished group-likes to obtain a formula for the square of the antipode, see \cite{Kauffman1993}.
	The next lemma follows from \cite[Proposition 6]{Kauffman1993}.

	\begin{lemma}
		Let $H$ be a pointed Hopf algebra, i.e. a Hopf algebra whose simple comodules are one-dimensional. If the dimension of $H$  is odd, it has a pair in involution $(l,\beta)$ such that $g \defeq l^{-2}$ and $\alpha \defeq \beta^2$ are the distinguished group-likes of $H^\circ$ and $H$, respectively.
	\end{lemma}
	
	Let us conclude this section with a remark on the representation theoretic viewpoint on pairs in involution. 
	Given a finite-dimensional Hopf algebra $H$ one can associate to it its category of finite-dimensional \YetterDrinfeld modules, see Section \ref{sec: Doubles}. 
	It is a \emph{rigid} category; i.e it is monoidal together with a notion of duality  compatible with its monoidal structure.
	In this context a pair in involution corresponds to, and can be reconstructed from, a \emph{monoidal} natural isomorphism between the identity functor and the functor which maps objects and morphisms to their biduals. 
	A  rigid category admitting such a structure is called \emph{pivotal}. In Sec\-tion \ref{sec: ADDandPiv} we will discuss the correspondence between pairs in involution and the pivotality of \YetterDrinfeld modules.	
	
	\section{Generalised Taft algebras} \label{sec: GTA}
	The strategy behind defining generalised Taft algebras is as follows. Fix a finite cyclic group $G$. Choose a \YetterDrinfeld module $V$ over the group algebra $\k \! G$ whose braiding is subject to certain relations.
	The bosonisation of the Nichols algebra of $V$  along $\k \! G$ yields another Hopf algebra. 
	This will be referred to as the coopposite of a generalised Taft algebra, see Definition \ref{def: GTA}. 
	A presentation in terms of generators and relations is obtained in Theorem \ref{thm: GTAGenerators}.

 	\subsection{\YetterDrinfeld and anti-\YetterDrinfeld modules}\label{sec: Doubles} Unless stated otherwise, every Hopf algebra in this section is assumed to have an invertible antipode.
 	
	The next Definition agrees with \cite[Definition 7.15.2]{Etingof2015}.
	\begin{definition}
		A \emph{\YetterDrinfeld module} over a Hopf algebra $H$ is a $\k$-vector space $M$ together with a  module structure $\lact\from H \otimes M \rightarrow M$ and a comodule structure  $\lcoact\from M \rightarrow H \otimes M$ satisfying the compatibility condition
		\begin{equation}
			\lcoact(h\lact m) = h_{(1)} m_{(-1)} S(h_{(3)}) \otimes h_{(2)} \lact m_{(0)}, \qquad \forall h\in H,m \in M.	\tag{YD}
		\end{equation} 
		A linear map $f\from M \rightarrow N$ between \YetterDrinfeld modules  is called a \emph{morphism of \YetterDrinfeld modules} if it is both a module and comodule morphism.
	\end{definition}

	\begin{remark}
		The \YetterDrinfeld modules over a Hopf algebra $H$ form the category $^H_H\mathcal {YD}$. The diagonal action and coaction of $H$ define a monoidal structure on it, see \cite[Chapter 7.15]{Etingof2015}.
		If the antipode of $H$ is invertible, 
		the natural isomorphism
		\begin{equation}
			\begin{aligned}
				&\sigma_{M,N} \from M \otimes N \rightarrow N \otimes M, 
				& m\otimes n \mapsto m_{(-1)}\lact n \otimes m_{(0)}. 
			\end{aligned}
		\end{equation}
		turns $_H^H\mathcal {YD}$ into a \emph{braided category}. This is explained in Chapters 7 and 8 of \cite[]{Etingof2015}.
	
		\YetterDrinfeld modules over a finite-dimensional Hopf algebra $H$ coincide with modules over its \emph{Drinfeld double} $D(H)$, see \cite[Chapter IX]{Kassel1995} \footnote{ \label{ftnote: Dfdouble}
			The definition of the Drinfeld and anti-Drinfeld double given here varies from \cite[Chapter IX.4]{Kassel1995} and \cite[Proposition 4.1]{Hajac2004} to accomodate our choice of (anti-)\YetterDrinfeld modules.
			In terms of the literature listed above our definition would read as
			$D(H^{\cop})^{\cop}$ and $A(H^{\cop})^{\cop}$. 
		}. 
		It is the vector space $H^\circ \otimes H$ whose Hopf algebra structure is for $g,h \in H$ and $\alpha, \beta \in H^\circ$ defined by
		\begin{equation}
			\begin{aligned}
				(\alpha \otimes g)(\beta \otimes h)  &
				\defeq \beta_{(1)}(S(g_{(1)})) \beta_{(3)} (g_{(3)}) \beta_{(2)} \alpha \otimes g_{(2)} h,
				\\
				\Delta(\alpha \otimes g)  &
				\defeq (\alpha_{(1)} \otimes g_{(1)}) \otimes (\alpha_{(2)} \otimes g_{(2)}),
				\\
				S(\alpha \otimes g )  &
				\defeq \alpha_{(1)}(g_{(1)}) \alpha_{(3)}(S(g_{(3)})) S^{-1}(\alpha_{(2)}) \otimes S(g_{(2)}).
			\end{aligned} \tag{DD}
		\end{equation}
	\end{remark}
	\begin{remark}
		The \emph{anti-Drinfeld double} $A(H)$ of a finite-dimensional Hopf algebra $H$ was introduced in \cite[Proposition 4.1]{Hajac2004}
		\textsuperscript{\ref{ftnote: Dfdouble}}. It is a 
		comodule algebra over the Drinfeld double $D(H)$.
		As a vector space it is $H^\circ\otimes H$.  The multiplication and coaction are given by
		\begin{equation}
			\begin{aligned}
			(\alpha \otimes g)(\beta \otimes h) & \defeq
			 \beta_{(1)}(S(g_{(1)})) \beta_{(3)} (S^{-2}(g_{(3)})) \beta_{(2)} \alpha \otimes g_{(2)} h, \\
			 \lcoact (\alpha \otimes g) & \defeq 
			 (\alpha_{(1)} \otimes g_{(1)}) \otimes (\alpha_{(2)} \otimes g_{(2)}),
			 \end{aligned}\tag{ADD}
		\end{equation}
		for $g,h \in H$ and $\alpha \in H^\circ$.
		Modules over $A(H)$ correspond to \emph{anti-\YetterDrinfeld modules}, see \cite[Definition 2.1]{Hajac2004}.
		That is, triples $(M,\lact, \lcoact)$ comprising a vector space $M$, an action $\lact \from H \otimes M \rightarrow M$ and a coaction $\lcoact \from M \rightarrow H \otimes M$ such that
		\begin{equation}\label{eq: AYD}
			\lcoact(h \lact m) = h_{(1)} m _{(-1)} S^{-1}(h_{(3)}) \otimes h_{(2)} \lact m_{(0)} , \qquad \forall h\in H,m \in M.
			\tag{AYD}
		\end{equation}	
		In general, anti-\YetterDrinfeld modules do not form a monoidal category but a module category over the \YetterDrinfeld modules \cite[Lemma 2.3]{Hajac2004}.
	\end{remark}
	Our interest in anti-\YetterDrinfeld modules is due to the following unpublished result by Hajac and Sommerh\"auser. We include a proof for the readers convenience.	
	\begin{theorem}\label{thm: PiiEquivIso}
		Suppose $H$ is a finite-dimensional Hopf algebra. The following are equivalent:
		\begin{enumerate}
			\item $A(H)$ and $D(H)$ are isomorphic as algebras.
			\item There exists a one-dimensional anti-\YetterDrinfeld module.
			\item $H$ has a pair in involution.
		\end{enumerate}
	\end{theorem}

	\begin{proof}
		Assume $f \from A(H) \rightarrow D(H)$ to be an isomorphism of algebras. 
		The ground field $\k$ considered as a module over $D(H)$ becomes a one-dimensional $A(H)$ module by pulling back the action along $f$. Hence, $(1)$ implies $(2)$. 
		
		Let $(\k ,\lact, \lcoact)$ be a one-dimensional anti-\YetterDrinfeld module. Its action and coaction are implemented by group-like elements $\beta^{-1} \in \Gr(H^\circ)$ and $l \in \Gr(H)$.
		We identify $H \cong H \otimes \k$ and observe
		\begin{equation} \label{eq: equivFormulaforPii}
			\beta^{-1}(h_{(2)})S(h_{(1)}) l  
			= S(h_{(1)}) \lcoact(h_{(2)} \lact 1) 
	 		\overset{ \text{\eqref{eq: AYD}}} { = } 
	 		\beta^{-1}(h_{(1)}) l S^{-1}(h_{(2)})
	 		\quad \forall h \in H.
		\end{equation}
		Applying $S$ to both sides proves \eqref{eq: equivFormulaforPii} to be equivalent to $(l, \beta)$ being a pair in involution.	Thus, $(3)$ follows from $(2)$.
		
		Given a pair in involution 	$(l,\beta)$ we define the linear map $f \from A(H) \rightarrow D(H)$,
		$\alpha \otimes g \mapsto \alpha_{(2)}(l) \beta^{-1}(g_{(2)}) \; \alpha_{(1)} \otimes g_{(1)}$ and compute for $(\alpha \otimes g), (\gamma \otimes h) \in A(H)$
		\begin{equation*}
			\begin{aligned}
				f( &(\alpha \otimes g) (\gamma \otimes h) ) =
				\gamma_{(1)}(S(g_{(1)})) 
				\gamma_{(3)}(S^{-2}(g_{(3)})) \;
				f(\gamma_{(2)} \alpha \otimes g_{(2)}h)
				\\
				& = \gamma_{(1)}(S(g_{(1)})) 
				\gamma_{(4)}(S^{-2}(g_{(4)}))
				\gamma_{(3)}(l)  
				\alpha_{(2)}(l) 
				\beta^{-1}(g_{(3)}) 
				\beta^{-1}(h_{(2)}) \;
				\gamma_{(2)}\alpha_{(1)} \otimes g_{(2)}h_{(1)}
				\\
				& = \gamma_{(1)}(S(g_{(1)})) 
				\gamma_{(3)}
				(\beta^{-1}(g_{(3)}) l S^{-2}( g_{(4)}))  
				\alpha_{(2)}(l)
				\beta^{-1}(h_{(2)}) \;
				\gamma_{(2)}\alpha_{(1)} \otimes g_{(2)}h_{(1)}
				\\
				\overset{(\ref{eq: AntipodeCondition})}&{=} 
				\gamma_{(1)}(S(g_{(1)}))
				\gamma_{(3)}(g_{(3)})  
				\gamma_{(4)}(l)
				\alpha_{(2)}(l)
				\beta^{-1}(g_{(4)}) 
				\beta^{-1}(h_{(2)})\;
				\gamma_{(2)}\alpha_{(1)} \otimes g_{(2)}h_{(1)}
				\\
				& = f( (\alpha \otimes g) ) f( (\gamma \otimes h) ).
			\end{aligned}
		\end{equation*}
		This proves $f$ to be a morphism of algebras.
		Its inverse is $f^{-1} \from D(H) \rightarrow A(H)$, 
		$\alpha \otimes g \mapsto \alpha_{(2)}(l^{-1}) \beta(g_{(2)}) \; \alpha_{(1)} \otimes g_{(1)}$. Therefore, $(3)$ implies $(1)$.
	\end{proof}	

	Our next result shows that the existence of a pair in involution corresponds to a suitably strong notion of Morita equivalence between the Drinfeld and anti-Drinfeld double.
	Given an algebra $A$ we write $\Forg_A \from A\text{-Bimod} \rightarrow A\text{-Mod}$ for the forgetful functor from the category of bimodules over $A$ to the category of left $A$ modules.
	
	\begin{lemma}\label{lem: MoritaDfAntiDf}
		Let $H$ be a finite-dimensional Hopf algebra. 
		It is equivalent:
		\begin{enumerate}
			\item $H$ has a pair in involution.
			\item There are $\k$-linear equivalences of categories 
			$F \from D(H)\text{\normalfont{-Mod}} \rightarrow A(H)\text{\normalfont{-Mod}}$ and
			$G \from D(H)\text{\normalfont{-Bimod}} \rightarrow A(H)\text{\normalfont{-Bimod}}$ such that a natural isomorphism $\eta \from F \; \Forg_{D(H)} \rightarrow \Forg_{A(H)} \; G$ exists.			
		\end{enumerate}
	\end{lemma}

	\begin{proof}
		Suppose $H$ has a pair in involution. 
		By Theorem \ref{thm: PiiEquivIso} there exists an isomorphism of algebras $f \from A(H) \rightarrow D(H)$. 
		Let $F \from D(H)\text{-Mod} \rightarrow A(H)\text{-Mod}$
	 	be the functor that identifies modules over $D(H)$ with modules over $A(H)$ by pulling back the action along $f$ and define	
	 	$G \from D(H)\text{-Bimod} \rightarrow A(H)\text{-Bimod}$ likewise. 
		Both, $F$ and $G$ are $\k$-linear equivalences of categories and $F \; \Forg_{D(H)} = \Forg_{A(H)} \; G$.
		
		Conversely, assume $F$, $G$ and $\eta$ to be as described above. 
		Let $X_\text{bi} \defeq (X,\lact, \ract)$ be a bimodule over $D(H)$.
		Set $X_\text{l} \defeq \Forg_{D(H)}(X_{\text{bi}})=(X,\lact)$ and write $Y \defeq F(X_\text{l})$.
		The module endomorphisms $\End{A(H)}{Y}$ themselves become a module over $A(H)$ via
		\begin{equation*}
			\tilde\lact \from A(H) \otimes \End{A(H)}{Y} \rightarrow \End{A(H)}{Y}, \quad
			(a \tilde{\lact} \phi) (x) \defeq
			\phi(\eta^{-1}_{X_\text{bi}}(\eta_{X_\text{bi}}^{\vphantom{-1}} (x) \ract a )).
		\end{equation*}
		As $F$ is a $\k$-linear equivalence of categories 
		$\End{D(H)}{X_\text{l}} \cong \End{A(H)}{Y}$ 
		as $\k$-vector spaces.
		Choose $X_{\text{bi}} \defeq {}_\epsilon \! \k _\epsilon$ to be the trivial bimodule over $D(H)$. 
		Then $X_\text{l} = {}_\epsilon \!\k$ is the trivial module over $D(H)$ and $\End{A(H)}{F({}_\epsilon \!\k)}$ is a one-dimensional module over the anti-Drinfeld double. The existence of a pair in involution follows from Theorem \ref{thm: PiiEquivIso}.
	\end{proof}

	\subsection{Nichols algebras and bosonisations}
	We follow the survey articles \cite{Andruskiewitsch2017, Andruskiewitsch2002} in recalling some aspects of Nichols algebras of diagonal type. 
	Until the end of this subsection we fix a Hopf algebra $H$ with invertible antipode.
	\newline
		
	The definition of Hopf algebras generalises naturally to braided monoidal categories. 
	A Hopf algebra $R$ in the category of \YetterDrinfeld modules over some Hopf algebra $H$ is referred to as a \emph{braided} Hopf algebra. Our next definition follows \cite[Definition 2.1]{Andruskiewitsch2002} almost verbatim. 
	\begin{definition}\label{def: NicholsAlgebra}
		Let $V$ be a \YetterDrinfeld module over $H$.
		A \emph{Nichols algebra} of $V$ is a braided graded Hopf algebra $\mathcal B(V)= \oplus_{n\geq 0} \mathcal B(V)(n)$ satisfying
		\begin{enumerate}
			\item $\mathcal B(V)(0)= \k$,
			\item $\mathcal B(V)(1) = \Pr(\mathcal B(V)) \cong V$ and
			\item $\mathcal B(V)$ is generated as an algebra by $B(V)(1)$. 
		\end{enumerate} 
	\end{definition}	
	Proposition 2.2 of \cite{Andruskiewitsch2002} asserts the existence of Nichols algebras and their uniqueness up to isomorphism. By abuse of notation we will speak of \emph{the} Nichols algebra in the following.  Aa explained in \cite[Definition 1.6]{Andruskiewitsch2002}, different types of Nichols algebras are distinguished in terms of their braidings.
	
	\begin{definition}\label{def: diagonalType}
		Let $V$ be a $\theta$-dimensional \YetterDrinfeld module over $H$.
		Write $\sigma \defeq \sigma_{V,V} \from V\otimes V \rightarrow V \otimes V$ for its braiding. 
		An ordered $\k$-basis $\{v_1,\dotsc, v_\theta \}$ of $V$ is said to be of \emph{diagonal type} if 
		\begin{equation}\label{eq: diagonalBraiding}
			\sigma (v_i \otimes v_j) = \mathfrak q_{ij} v_j \otimes v_i, \qquad \mathfrak q_{ij}\in \k \text{ for } 1 \leq i,j\leq \theta.
		\end{equation}
		Accordingly, $V$ and $\mathcal B(V)$ are referred to be of diagonal type if $V$ has an ordered basis of diagonal type.
		The matrix 
		$(\mathfrak q_{ij})\in \k^{\theta \times \theta}$
		is called the \emph{matrix of the braiding}.
	\end{definition}

	Finite-dimensional Nichols algebras of diagonal type were classified by Heckenberger in terms of generalised Dynkin diagrams, see \cite{Heckenberger2009}. 
	The Nichols algebra part of a generalised Taft algebra corresponds to a  diagram of $A_1\times A_1 \cong D_2$-type. 
	
	\begin{definition}
		A \YetterDrinfeld module $V$ over $H$ is of \emph{$D_2$-type} if
		a basis of diagonal type exists whose
		matrix of the braiding $(\mathfrak q_{ij})\in \k^{2\times 2}$
		satisfies:
		\begin{equation}\label{eq: D2Conditions}
		\text{Its entries are roots of unity,
			$\mathfrak q_{11}, \mathfrak q_{22} \neq 1 $ and 
			$\mathfrak q_{12} \mathfrak q_{21} = 1$}. 
		\end{equation}
		Likewise, its Nichols algebra $\mathcal B(V)$ is also referred to as of $D_2$-type.
	\end{definition}
	
	Nichols algebras $\mathcal B(V)$ of $D_2$-type are also referred to as \emph{quantum planes}, see \cite{Andruskiewitsch1998a}.

	\begin{remark}\label{rem: IndOfBasis}
		A direct computation shows that  whether a \YetterDrinfeld module is of $D_2$-type does not depend on the choice of basis of diagonal type.
	\end{remark}

	The bosonisation or biproduct of a braided Hopf algebra $R$ over $H$ equips the vector space $R\otimes H$ with the structure of a Hopf algebra, see \cite[Theorems 11.5.7 and 11.6.9]{Radford2012}.
	To distinguish between the coaction and comultiplication of $R$ we use a slight variation of Sweedler notation and write $r^{(1)} \otimes r^{(2)} \defeq \Delta (r)$ for $r\in R$.

	\begin{definition}\label{def: bos}
		Let $H$ be a Hopf algebra whose antipode is invertible and $R$ a braided Hopf algebra in $^H_H \mathcal{YD}$. 
		The \emph{bosonisation} of $R$ by $H$ is the Hopf algebra $R\#H$, whose underlying vector space is $R \otimes H$ and whose multiplication, comultiplication and antipode is defined for $g,h \in H$ and $r,s \in R$ by
		\begin{equation}
			\begin{aligned}
				 (r\otimes g)(s\otimes h) &\defeq 
				 r(g_{(1)}\lact s) \otimes g_{(2)}h,
				\\
				\Delta(r\otimes g) &\defeq 
				r^{(1)}\otimes (r^{(2)})_{(-1)}g_{(1)} \otimes (r^{(2)})_{(0)}\otimes g_{(2)},
				\\
				S(r\otimes g) &\defeq 
				S_H(r_{(-2)} g_{(2)})\lact S_R(r_{(0)})\otimes S_H(r_{(-1)} g_{(1)}).
			\end{aligned}
		\end{equation}		
	\end{definition}

	\subsection{Generalised Taft algebras} 
	We define now the main object under investigation, generalised Taft algebras. In 
	Theorem \ref{thm: GTAGenerators} we obtain a presentation in terms of generators and relations.

	\begin{definition}\label{def: GTA}
		Let $V$ be a  \YetterDrinfeld module of $D_2$-type over the group algebra $\k \! G$ of a finite cyclic group $G$.
		We call $(\mathcal B (V) \#\k \! G)^{\cop}$
		a \emph{generalised Taft algebra}.
	\end{definition}

	In the above definition the coopposite is chosen to match the definition of Taft algebras as given for example in \cite[p. 113]{Kauffman1993}.
	
	Let $N \geq 2$ be a natural number. We write $\mathbb Z_N \defeq \; \factor{\mathbb Z}{N\mathbb Z}$.
	\begin{definition}
		A matrix 
		$(a_{ij})
		\in \mathbb Z_N^{2 \times 2}$
		whose entries satisfy  
		modulo $N$
		\begin{equation}\label{eq: GTAcondition}
			a_{11} a_{12} \neq 0, \qquad
			a_{21} a_{22} \neq 0, \qquad
			a_{11} a_{22} + a_{12} a_{21} = 0 
		\end{equation} 
		is called 
		a \emph{parameter matrix} of a generalised Taft algebra.
	\end{definition}
	
	Parameter matrices allow us to define  \YetterDrinfeld modules of $D_2$-type.
	
	\begin{lemma}
		Let $N\geq 2$ and $(a_{ij})
		\in \mathbb Z_N^{2 \times 2}$. Fix a generator $g \in \mathbb Z_N$ and a primitive $N$-th root of unity $q\in \k$.
		The \YetterDrinfeld module $V \defeq \spanset_{\k}\{x, y\}$, defined by 
		\begin{gather} \label{eq: actAndCoact}
			\lcoact(x) \defeq g^{a_{11}} x, \qquad
			\lcoact(y) \defeq g^{a_{21}} y, \qquad
			g \lact x \defeq q^{a_{12}} x, \qquad
			g \lact y \defeq q^{a_{22}} y,	
		\end{gather}
		generates a generalised Taft algebra $(\mathcal B(V) \# \k \! \mathbb Z_N)^{\cop}$ if and only if $(a_{ij})$ is a parameter matrix.
	\end{lemma}
	\begin{proof}
		The matrix of the braiding of $V$ is, with respect to the ordered basis $\{x, y\}$, given by
		$\mathfrak q_{ij} \defeq q^{a_{j1}a_{i2}}$. Identifying $q$ with a generator of $\mathbb Z_N$ shows that $V$ is of $D_2$-type if and only if $(a_{ij})$ is a parameter matrix.
	\end{proof}

	Given a matrix $(a_{ij})\in \mathbb Z_N^{2 \times 2}$, we write  
	$N_x \defeq  \ord(a_{11}a_{12})$ and $N_y \defeq \ord(a_{21} a_{22})$ for the orders of $a_{11}a_{12}$ and $a_{21}a_{22}$ in $\mathbb Z_N$, respectively. 
	
	\begin{theorem}\label{thm: GTAGenerators}
		Let $H$ be a generalised Taft algebra. 
		Then there exists an integer $N \geq 2$, a parameter matrix 
		$(a_{ij})
		\in \mathbb Z_N^{2 \times 2}$ and a primitive $N$-th root of unity $q\in \k$ 
		such that $H$ is generated by elements
		 $g, x, y \in H$ subject to the relations
		 \begin{subequations}
		\begin{gather}
			\begin{gathered} \label{eq: NichAlgPartRelations}				
				x^{N_x} = 0, \qquad
				y^{N_y} = 0, \qquad
				xy = q^{a_{11}a_{22}} yx,
			\end{gathered}
				\\
			\begin{gathered} \label{eq: GrpAlgPartRelations}
				g^N = 1, \qquad
				gx = q^{a_{12}}xg, \qquad
				gy = q^{a_{22}}yg,
			\end{gathered}
			\\
			\begin{gathered}\label{eq: CoalgPartRelations}
				\Delta(g) = g \otimes g, \quad
				\Delta(x) = 1 \otimes x + x \otimes g^{a_{11}},  \quad
				\Delta(y) = 1 \otimes y + y \otimes g^{a_{21}},
			\end{gathered}
			\\
			\begin{gathered}\label{eq: AntipodeRelations}
				S(g) = g^{-1}, \qquad
				S(x) = - xg^{-a_{11}}, \qquad
				S(y) = - yg^{-a_{21}}.
			\end{gathered}
		\end{gather}
		\end{subequations}
	\end{theorem}
	\begin{proof}
		By definition $H = (\mathcal B (V) \# \k \!\mathbb Z_N)^{\cop}$ for some $ N\geq 2$.
		We fix a generator $\\overline g\in \mathbb Z_N$ and a primitive $N$-th root of unity $q \in \k$. 
		 As a vector space $V$ admits
		an ordered basis $\{\overline x, \overline y\}$ such that
		a  matrix $(a_{ij})
		\in \mathbb Z_N^{2\times 2}$ exists which implements the action and coaction as in Equation \eqref{eq: actAndCoact}, see \cite[Remark 1.5]{Andruskiewitsch2002}.
		The preceding lemma shows that $(a_{ij})$ needs to be a parameter matrix.
		Every element $z \otimes h \in H$ can be factorised into the product
		$(z \otimes 1)(1 \otimes h)$, for $z \in \mathcal B(V)$ and $h \in \k\! \mathbb Z_N$, implying that $g \defeq 1 \otimes \overline g$, $x \defeq \overline x \otimes 1$ and $y \defeq  \overline y \otimes 1$ generate $H$.
		The relations \eqref{eq: NichAlgPartRelations} follow from \cite[Corollary 8.1]{Heckenberger2007}.
		The definition of the multiplication of a bosonisation and $\overline g^N = 1$  imply the relations \eqref{eq: GrpAlgPartRelations}. 
		The coproduct and antipode of the generators are obtained by the respective formulas in Definition \ref{def: bos}.
	\end{proof}

	\begin{convention}
		We fix some of the notation of this section. 
		From now onwards $N\geq 2$ denotes an integer and $1 \neq q\in \k$ a primitive $N$-th root of unity. 
		We write $A \defeq (a_{ij}) \in \mathbb Z_N^{2 \times 2}$ for a parameter matrix and $H \defeq H_q(A)$ for its associated generalised Taft algebra, generated by elements $g$, $x$ and $y$. 
		In particular, $g$ generates the group $\Gr(H)$ and $x$ and $y$ are nilpotent of degree $N_x$ and $N_y$, respectively.
	\end{convention}
	
	\section{Hopf algebras without pairs in involution}\label{sec: HAWithoutPII}
	
	We show that similar to Taft algebras their generalisations form a class of  basic, pointed Hopf algebras, which is closed under duality.
	Our main result, Theorem \ref{thm: Classification}, states necessary and sufficient conditions for these Hopf algebras
	to omit  pairs in involution.
	Thereafter, we investigate how various properties of generalised Taft algebras affect the existence of such pairs. 
	In Lemma \ref{lem: Examples} we construct an infinite family of examples of finite-dimensional Hopf algebras without pairs in involution. 
	\newline
		
	\subsection{Properties of generalised Taft algebras}
			
	To prove that the class of generalised Taft algebras is closed under duality, we identify generators of the dual. 
	For a generalised Taft algebra $H\defeq H_q(A)$,  generated by $g, x, y \in H$, we set $\xi , \psi, \phi\from H \rightarrow \k$,
	\begin{equation}\label{eq: dualGenerators}
		\begin{gathered}
			\xi(x^i y^j g^l) \defeq q^{-l} \delta_{i=j=0}, 
			\quad 
			\psi(x^i y^j g^l) \defeq q^{-a_{12} l} \delta_{i=1, j=0}, 
			\\	
			\phi(x^i y^j g^l) \defeq q^{-a_{22} l} \delta_{i=0 , j=1}.
		\end{gathered}
	\end{equation}
	Given a matrix $A \in \mathbb Z_N^{2\times 2}$ we write $A_{21}$ for the matrix obtained by interchanging the first with the second column of $A$. The next result is implied by {\cite[Proposition 3.1]{Nenciu2004}}.
	\begin{theorem}\label{thm: Duality}
		Let $H \defeq H_q{(A)}$ be 
		a generalised Taft algebra. Write $\bar g$, $\bar x$ and $\bar y$ for the generators of $H_q{(A_{21})}$.
		The map $\Theta \from H_q{(A_{21})} \rightarrow H^\circ $,
		$\Theta(\bar g) \defeq \xi$, $\Theta(\bar x) \defeq \psi$ and $\Theta(\bar y) \defeq \phi$ is an isomorphism of Hopf algebras.
	\end{theorem}

	A Hopf algebra is called \emph{pointed} if every simple comodule is one-dimensional; if every simple module is one-dimensional it is called \emph{basic}.
	
	\begin{lemma}
		Generalised Taft algebras are pointed and basic.
	\end{lemma}

	\begin{proof}
		The first claim follows from \cite[Proposition 4.4.9]{Radford2012}. The second is a consequence of the former one and Theorem \ref{thm: Duality}. 
	\end{proof}

	We determine the left integrals and distinguished group-likes of  generalised Taft algebras. The latter will prove useful in the study of the square of the antipode. This is a standard exercise, see for example \cite[Section 2.12]{Andruskiewitsch2017}.

	\begin{lemma}\label{lem: integrals}
		Let 
		$H \defeq H_q(A)$ be a generalised Taft algebra. The left integrals of $H$ and $H^\circ$ are up to scalar multiplication
		\begin{equation*}
			\Lambda \defeq \left(\sum\nolimits_{i=0}^{N-1} g^i \right)x^{N_x-1}y^{N_y-1} \in H \text{ and }
			\Upsilon \defeq \left(\sum\nolimits_{i=0}^{N-1} \xi^i \right) \psi^{N_x-1} \phi^{N_y-1} \in H^{\circ}.
		\end{equation*} 
		The elements $\xi^{-(a_{12} +a_{22})} \in H^\circ$ and $g^{-(a_{11} +a_{21})} \in H$ are the
		distinguished group-likes of $H$ and $H^{\circ}$, respectively.
	\end{lemma}

	\begin{proof}
		By multiplying $\Lambda$ with the generators of $H$ we see that it is a left integral and that $\xi^{a_{12} (N_x-1) + a_{22} (N_y-1)}$ is the distinguished group-like of $H$. 
		Modulo $N$ we have $a_{12} N_x = a_{22} N_y =0$ and therefore
		$\xi^{a_{12} (N_x-1) + a_{22} (N_y-1)} = \xi^{-(a_{12} +a_{22})}$. 
		The results for $H^\circ$ follow by applying the isomorphism of Theorem \ref{thm: Duality}.
	\end{proof}
	
	The motivation behind Kauffman's and Radford's study of the square of the antipode, see \cite{Kauffman1993}, was understanding Hopf algebras which give rise to knot invariants. 
	This necessarily requires the Hopf algebra to be \emph{quasitriangular}. 
	That is, roughly speaking, an encoding of the notion of braidings on the level of Hopf algebras, see \cite[Chapters VIII and XIII]{Kassel1995}. 
	The next Lemma is implied by Theorem 3.4 of \cite[]{Nenciu2004}.
	\begin{lemma}\label{lem: quasitriangular}
		Let $H \defeq H_q(A)$ be a generalised Taft algebra with generators $g$, $x$, $y$. Recall that $N_x, N_y \in \mathbb N$ were the minimal positive integers such that $x^{N_x} = 0$ and $y^{N_y} = 0$. 
		Then $H$ is quasitriangular if and only if $N$ is divisible by $2$,  $N_x = N_y = 2$ and $a_{11} = a_{21} = \frac{N}{2}$.	
	\end{lemma}	
	
	\subsection{Pairs in involution as solutions of Diophantine equations}
	
	To find necessary and sufficient criteria for generalised Taft algebras to not admit a pair in involution we study the behaviour of the square of the antipode.
	Fix a generalised Taft algebra $H \defeq H_q(A)$.
	Its square of the antipode is determined by
	\begin{equation*}
		S^2(g) = g, \qquad
		S^2(x) = q^{a_{11}a_{12}} x, \qquad
		S^2(y) = q^{a_{21}a_{22}} y.
	\end{equation*}
	Likewise the fourth power of the antipode is
	\begin{equation*}
		S^4(g) = g, \qquad
		S^4(x) = q^{2a_{11}a_{12}} x, \qquad
		S^4(y) = q^{2a_{21}a_{22}} y.
	\end{equation*}	
	Now, consider the family of Hopf algebra automorphisms $T_{(l,\beta)} \from H \rightarrow H$ which is indexed by a pair of group-likes $l \in \Gr(H)$, $\beta \in \Gr(H^\circ)$ and defined via 
	\begin{equation*}
		T_{(l,\beta)}(h)= \beta (h_{(3)}) \beta^{-1} (h_{(1)}) \;  l h_{(2)} l^{-1},
		\qquad \text{ for all } h \in H.
	\end{equation*}
	Definition \ref{def: PII} states that $H$ has a pair in involution if and only if there exists a pair 
	of group-like elements 
	$(l,\beta)$ such that $T_{(l,\beta)} = S^2$. 
	Every group-like $l\in \Gr(H)$ and character $\beta \in \Gr(H^\circ)$ can uniquely be written as  $l=g^d$ and 
	$\beta= \xi^{-c}$ with $c,d \in \mathbb Z_N$. 
	Evaluating $T_{(g^d, \xi^{-c})}$ on the generators yields
	\begin{equation*}
		T_{(g^d, \xi^{-c})}(g) = g, \qquad
		T_{(g^d, \xi^{-c})}(x) = q^{a_{11}c +a_{12}d} x, \qquad
		T_{(g^d, \xi^{-c})}(y) = q^{a_{21}c +a_{22}d} y.
	\end{equation*}
	 Identifying the $N$-th roots of unity with $\mathbb Z_N$ via
	 $q \mapsto 1$ implies our next theorem. 
	 
	 \begin{theorem}\label{thm: ModPairAsEquation}
	 	Let $q $ be a primitive $N$-th root of unity. 
	 	A generalised Taft algebra $H_q(A)$ has a pair in involution if and only if
	 	$c,d \in \mathbb Z_N$ exist such that modulo $N$
	 	\begin{equation}\label{eq: ModPair}
		 	a_{11}c + a_{12}d = a_{11}a_{12}, \qquad 
		 	a_{21}c + a_{22}d = a_{21}a_{22}.
	 	\end{equation}
	 	The pair is modular if and only if it additionally satisfies modulo $N$
	 	\begin{equation}
	 		cd = 0.
		\end{equation}
	\end{theorem}

 	\begin{remark}\label{rem: SolutioninOddCase}
 		Lemma \ref{lem: integrals} and the discussion prior to the above theorem imply that for a given parameter matrix $(a_{ij}) \in \mathbb Z_N^{2 \times 2}$ of a generalised Taft algebra  the integers $c'\defeq a_{12} + a_{22}$, and  $d'\defeq a_{11} +a_{21}$ satisfy modulo $N$
 		\begin{equation*}
	 		a_{11}c' + a_{12}d' = 2 a_{11}a_{12}, \qquad 
	 		a_{21}c' + a_{22}d' = 2 a_{21}a_{22}.
 		\end{equation*}
 		If $N$ is odd, $2$ is invertible modulo $N$. In this case 
 		\eqref{eq: ModPair} has a solution.
 	\end{remark} 
 
 	Consequently, a generalised Taft algebra without pairs in involution necessarily needs to have a group of group-likes of even order $N$. 
 	In fact, the existence of such a pair depends on the behaviour of the entries of the parameter matrix modulo $2^n$, where $n\in \mathbb N_0$ is the maximal integer such that $2^n$ divides $N$.

 	\begin{definition}\label{def: PowersAndCoefficients}
 		Suppose $2^{n}\cdot j \geq 2$, with $j$ odd, and let
 		$(a_{ij})\in \mathbb Z_{2^n j}^{2 \times 2}$
 		be a parameter matrix of a generalised Taft algebra.
 		The \emph{matrix of powers} associated to $(a_{ij})$ is the matrix 
 		$(\mathfrak a_{ij}) \in \mathbb N_0^{2 \times 2}$ whose entries are the minimal non-negative integers satisfying 
 		\begin{equation}\label{eq: coefAndPowerMatrix} 
 			2^{\mathfrak a_{ij}}\mu_{ij} = a_{ij} \mod 2^n, 
 			\qquad \text{ for } \mu_{ij}\in \mathbb Z_{2^n} 
 			\text { invertible or zero.}
 		\end{equation}
 		The matrix $(\mu_{ij})\in \mathbb Z_{2^n}^{2 \times 2}$ 
 		is called the \emph{matrix of coefficient} of $(a_{ij})$. 
 		The \emph{power of the coefficient matrix} is the minimal integer $\tau \in \mathbb N_0$ such that $2^\tau \nu= \det (\mu_{ij})$ modulo $2^n$, with $\nu\in \mathbb Z_{2^n}$ invertible or zero.	
	\end{definition}

	 \begin{theorem}\label{thm: Classification}
		Suppose $2^n \cdot j \geq 2$, with $j$ odd. 
		Let $q\in \k$ be a primitive $2^n j$-th root of unity and 
		$(a_{ij})$ a parameter matrix whose matrices of powers and coefficients are $(\mathfrak a_{ij})$ and $(\mu_{ij})$, respectively. Write $\tau$ for the power of $(\mu_{ij})$.
		It is equivalent: 
		\begin{enumerate}
			\item $H_q(a_{ij})$ has no pair in involution,
			\item 
			$n\geq 1$, 
			$\mu_{ij} \neq 0$ 
			for $1 \leq i, j \leq 2$, 
			$\mathfrak a_{11} + \mathfrak a_{22} < n$, $\mathfrak a_{11} \neq \mathfrak a_{21}$ and either $\tau > \min\{
			\mathfrak{a}_{ij} \mid 1 \leq i,j \leq 2\}$ or $\det (\mu_{ij}) = 0$ modulo $2^n$.
		\end{enumerate}
	\end{theorem}	
	\begin{proof}
		We prove $H$ having a pair in involution equivalent to at least one of the above conditions not being met.
		By Theorem \ref{thm: ModPairAsEquation} this amounts in solving 
		\begin{equation}\label{eq: PiiMatrix}
			\left(\begin{array}{cc|c} 
				a_{11} & a_{12} & a_{11} a_{12} \\
				a_{21} & a_{22} & a_{21} a_{22}
			\end{array}\right)
		\end{equation}
		modulo $2^n j$.
		As $\chi \from \mathbb Z_{2^n j} \rightarrow \mathbb Z_j \times \mathbb Z_{2^n}$, $(x \mod 2^n \cdot j) \mapsto (x \mod j, x \mod 2^n)$ is an isomorphism of rings
		this is equivalent to the solvability of the above equation 
		modulo $j$ and modulo $2^n$, respectively.
		In Remark \ref{rem: SolutioninOddCase}
		a solution modulo $j$ was given.
		Thus, $H$ has a pair in involution if and only if Equation \eqref{eq: PiiMatrix} is solvable modulo $2^n$.
		If $n=0$, this is trivially the case. Hence, we assume $n\geq 1$. 	
		Using the matrices of powers and coefficients, Equation \eqref{eq: PiiMatrix} can be expressed modulo $2^n$ as
		\begin{equation}\label{eq: Mod2powerUnnormalized}
			\left(\begin{array}{cc|c}  
				\mu_{11} 2^{\mathfrak a_{11}} & 
				\mu_{12} 2^{\mathfrak a_{12}} &
				\mu_{11}\mu_{12} 2^{\mathfrak a_{11}+ \mathfrak a_{12}} \\  
				\mu_{21} 2^{\mathfrak a_{21}} & 
				\mu_{22} 2^{\mathfrak a_{22}} &
				\mu_{21}\mu_{22} 2^{\mathfrak a_{21}+ \mathfrak a_{22}}
			\end{array}\right).
		\end{equation}	
		We may assume without loss of generality that
		$\mathfrak a_{11} \leq \mathfrak a_{ij}$ for $1 \leq i,j \leq 2$. 
			 
		If one of the terms on the right hand side is zero, a solution can be written down directly. For example, if $\mu_{12} = 0$ or 
		$2^{\mathfrak a_{11} + \mathfrak a_{12}} = 0$, take 
		$(0, \mu_{21} 2^{\mathfrak a_{21}})^{\T} \in \mathbb Z_{2^n}^2$.
		We therefore assume $\mu_{ij} \neq 0$ in the following.
		A less obvious implication follows from the identity
		$a_{11} a_{22} + a_{21} a_{12} = 0$ modulo $N$,
		which reads in our setting as
		\begin{equation}\label{eq: a1b2 + a2b1}
			\mu_{11} \mu_{22} 2^{\mathfrak a_{11}+ \mathfrak a_{22}} + 
			\mu_{12} \mu_{21} 2^{\mathfrak a_{12}+ \mathfrak a_{21}}  = 0 \mod 2^n.
		\end{equation}
		If $\mathfrak a_{11} + \mathfrak a_{22} \geq n$, or equivalently $\mathfrak a_{21} + \mathfrak a_{12} \geq n$, the second entry of the right hand side of \eqref{eq: Mod2powerUnnormalized} is zero. 
		Hence, we furthermore add $\mathfrak a_{11} + \mathfrak a_{22} < n$ to our list of assumptions, which in particular implies
		\begin{equation}\label{eq: mfrak a1b2 + a2b1}
			\mathfrak a_{11} + \mathfrak a_{22} = \mathfrak a_{21} + \mathfrak a_{12}.
		\end{equation}

		We transform \eqref{eq: Mod2powerUnnormalized} into upper triangular form by multiplying the second row with $\mu_{11}$ and then subtracting the first row $\mu_{21}2^{\mathfrak a_{21}- \mathfrak a_{11}}$-times from it.
		Using Equations \eqref{eq: a1b2 + a2b1} and \eqref{eq: mfrak a1b2 + a2b1}, this simplifies to
		\begin{equation}\label{eq: Mod2powerNormalizedBlockDiag}
		 	\left(\begin{array}{cc|c}  
		 		\mu_{11}2^{\mathfrak a_{11}} & 
		 		\mu_{12} 2^{\mathfrak a_{12}} &
		 		\mu_{11}\mu_{12} 2^{\mathfrak a_{11}+ \mathfrak a_{12}} \\  
		 		0 & 
		 		\det(\mu_{ij}) 2^{\mathfrak a_{22}} &
		 		(\mu_{22}\mu_{11})
		 		( \mu_{21} 2^{\mathfrak a_{21}- \mathfrak a_{11}} + \mu_{11})
		 		2^{\mathfrak a_{11}+ \mathfrak a_{22}}
		 	\end{array}\right).
		\end{equation}
		By the minimality of $\mathfrak{a}_{11}$ Equations \eqref{eq: Mod2powerUnnormalized} and \eqref{eq: Mod2powerNormalizedBlockDiag} have the same set of solutions and the solvability of \eqref{eq: Mod2powerNormalizedBlockDiag}  depends only on the existence of a solution for the equation displayed in its second row.
		We want to simplify this equation further.
		First, by the definition of the power of the coefficient matrix there exists an element $\nu \in \mathbb Z_{2^n}$, which is invertible or zero, such that
		$\det (\mu_{ij}) = 2^\tau \nu$. Second, we write $\rho \in \mathbb N$ for the minimal number such that 
		$ (\mu_{22}\mu_{11})
		( \mu_{21} 2^{\mathfrak a_{21}- \mathfrak a_{11}} + \mu_{11}) = 
		2^\rho \varpi$ for $\varpi \in \mathbb Z_{2^n}$, with $\varpi$ invertible or zero. 
		We divide the second row of Equation \eqref{eq: Mod2powerNormalizedBlockDiag} by $2^{\mathfrak a_{22}}$ and observe that it
		is solvable if and only if a $d \in \mathbb Z_{2^{n - \mathfrak a_{22}}}$ exists such that 
		\begin{equation}\label{eq: Mod2powerLowerBlock}
			2^\tau \nu \cdot  d = 2^{\mathfrak a_{11}+\rho} \varpi \mod 2^{n - \mathfrak a_{22}}.
		\end{equation}
		There are three cases which might occur.
		The right hand side might be zero and the equation is trivially solvable.
		This is the case if and only if $\varpi = 0$ or $\rho + \mathfrak a_{11} \geq n - \mathfrak a_{22}$. 
		One verifies that necessarily $\mathfrak a_{21} = \mathfrak a_{11}$ needs to hold.
		
		Now, let us assume that the right hand side is not zero but $\mathfrak a_{11} = \mathfrak a_{21}$.
		In particular, we have $\rho \geq 1$.
		By Equation \eqref{eq: a1b2 + a2b1} $\mu_{12} \mu_{21} = - \mu_{11} \mu_{22}$ modulo $2^{n -(\mathfrak a_{11} + \mathfrak a_{22})}$  and we can write
		$2^\tau \nu = \det(\mu_{ij})= 2(\mu_{11} \mu_{22} + \lambda 2^{n -(\mathfrak a_{11} + \mathfrak a_{22}+1)})$ for a $\lambda \in \mathbb Z$.
		Since the right hand side is non-zero we have $\mathfrak a_{11}+ \mathfrak a_{22} < n - \rho \leq n-1$ and $(\mu_{11} \mu_{22} + \lambda 2^{n -(\mathfrak a_{11} + \mathfrak a_{22} +1)})$ is, as the sum of an odd and even integer, invertible modulo $2^{n- \mathfrak a_{22}}$. 
		In other words, $\nu \neq 1$ and $1 =\tau \leq \rho$ and Equation \eqref{eq: Mod2powerLowerBlock} is solvable.
					
		Finally, consider the case where the right hand side is not zero and $\mathfrak a_{11} \neq \mathfrak a_{21}$.
		Then $\rho = 0$ and $\varpi \neq 0$. A solution exists if and only if 
		$\tau \leq \mathfrak a_{11}$ and $\nu \neq 0$.
	\end{proof}

	We study the effect of the order of the group of group-likes of a generalised Taft algebra on the existence of pairs in involution.
	
	\begin{lemma}
		Let $q\in \k$ be a primitive $2^n j$-th root of unity, with $j$ odd, and $H_q(A)$ a generalised Taft algebra without a pair in involution. Then $n \geq 2$ and $2^nj \neq 4$.
	\end{lemma}
	\begin{proof}
		Assume  $H_q(A)$ to not admit a pair in involution, which implies $n \geq 1$.
		Let $(\mathfrak a_{ij})$ be the matrix of powers associated to $A$.
		In case $n=1$, it is the zero matrix and by Theorem \ref{thm: Classification} a pair in involution exists.
		If $2^n j = 4$, note that we have  $0 = \mathfrak a_{11} \neq \mathfrak a_{21} = 1$. 
		The discussion prior to Equation \eqref{eq: mfrak a1b2 + a2b1} shows that
		$\mathfrak a_{11}+ \mathfrak a_{22} = \mathfrak a_{12} + \mathfrak a_{21}$ and, therefore, $\mathfrak a_{22} = 1$. But $\mathfrak a_{21} + \mathfrak a_{22} = 2$ implies $a_{12} a_{22} = 0$ modulo $4$, contradicting that $A$ is a parameter matrix.
	\end{proof}
	
	Conversely, we obtain for every natural number $4\cdot N \geq 8$ an example of a generalised Taft algebra without such a pair.
	
	\begin{lemma}\label{lem: Examples}
		Let $N \geq 2$ be a natural number and $q\in \k$ a primitive $4N$-th root of unity. The generalised Taft algebra $H_q
		\left(\begin{smallmatrix} 
					1 & 1 \\ 
					2 & -2 
		 \end{smallmatrix} \right)$ has no pair in involution.
	\end{lemma}

	\begin{proof}
		The matrices of powers and coefficients associated to $\left(\begin{smallmatrix} 
			1 & 1 \\ 
			2 & -2 
		\end{smallmatrix} \right)$ are 
		$(\mathfrak a_{ij}) =\left(\begin{smallmatrix}
			0 & 0 \\
			1 & 1
		\end{smallmatrix}\right)$ and 
		$(\mu_{ij}) = \left(\begin{smallmatrix}
			1 & 1 \\
			1 & -1
		\end{smallmatrix}\right)$. 
		Since $\det(\mu_{ij}) = -2$, Theorem \ref{thm: Classification} implies that $H_q\left(\begin{smallmatrix} 
			1 & 1 \\ 
			2 & -2 
		\end{smallmatrix} \right)$ does not admit a pair in involution.
	\end{proof}

	\begin{remark}
		Suppose $q\in \k$ to be a primitive $p$-th root of unity for $p>2$ a prime number. For any $s\in \mathbb Z_p \setminus \{0 \}$, the generalised Taft algebra 
		$H_q\left(\begin{smallmatrix} 
			1 & 1 \\ 
			s & -s 
		\end{smallmatrix} \right)$ has a pair in involution.
		In \cite{Halbig2019}, Kr\"ahmer and the author showed, however, that the modularity condition can only be satisfied if $s \in \{1, p-1\}$.
	\end{remark}
	
	Requiring the existence of additional structures on generalised Taft algebras might restrict the possible choices of parameter matrices severely. 
	We conclude this section by investigating two such cases and their influence on pairs in involution. 
	
	The next theorem, whose proof resembles \cite[Theorem 3]{Kauffman1993}, 
	will help us greatly in that respective. It
	links the existence of such pairs to their existence for the dual and Drinfeld double.
	\begin{theorem}\label{thm: equivPiiHDualDouble}
		Let $H$ be a finite-dimensional Hopf algebra. 
		If any of the Hopf algebras $H$, $H^\circ$ or $D(H)$ admit a pair in involution, then all do.
	\end{theorem}
	
	\begin{proof}
		Let $(l, \beta)$ be a pair in involution for the finite-dimensional Hopf algebra $H$.
		For any $\omega \in H^\circ$ and $h\in H$ we compute
		\begin{align*}
			S^2(\omega)(h) & = \omega (S^2(h)) 
			= \beta(h_{(3)}) \beta^{-1}(h_{(1)}) \omega (lh_{(2)} l^{-1}) \\ &
			=  \omega_{(1)}(l) \omega_{(3)}(l^{-1})\; (\beta^{-1} \omega_{(2)} \beta)(h).
		\end{align*}
		Thus, $(\beta^{-1}, l^{-1})$ is a pair in involution for $H^\circ$.
		The converse statement follows by the same argument and the fact that $H$ is isomorphic to its bidual.	
		
		For any element $\omega \otimes h \in D(H)$ we have
		\begin{equation} \label{eq: PiiDDimpliesPiv}
			\begin{aligned}
			S^2& (\omega \otimes h) 
			= (S^{-2}(\omega)\otimes 1)
			(\epsilon \otimes S^2(h)) \\ &
			= \omega_{(1)}(l^{-1}) \omega_{(3)}(l)
			\beta^{-1}(h_{(1)}) \beta(h_{(3)})
			( \beta \omega_{(2)} \beta^{-1} \otimes  lh_{(2)}l^{-1}) \\ &
			= \omega_{(1)}(l^{-1}) \omega_{(3)}(l) (\omega_{(2)}\beta^{-1} \otimes lh)(\beta \otimes l^{-1}) 
			= (\beta^{-1} \otimes l) (\omega \otimes h) (\beta \otimes l^{-1}). 
			\end{aligned}
		\end{equation}
		In other words, $(\beta^{-1}\otimes l, \epsilon_{D(H)})$ is a pair in involution for $D(H)$.
		
		Assume conversely that $\gamma \in \Gr(D(H)^\circ)$ and $c\in \Gr(D(H))$ constitute a pair in involution for $D(H)$.
		By \cite[Proposition 13.2.2]{Radford2012}, we have $c = \beta \otimes l$ for $\beta \in \Gr(H^\circ)$ and $l \in \Gr(H)$.
		As $\iota \from H \rightarrow D(H)$, $h \mapsto \epsilon \otimes h$ is an inclusion of Hopf algebras we obtain a character $\tilde \gamma \defeq \gamma \circ \iota \in \Gr(H^\circ)$.
		For any $h\in H$
		\begin{align*}
			S^2(\iota(h)) & = 
			\tilde \gamma( h_{(3)})
			\tilde \gamma^{-1}(h_{(1)})
			(\beta\otimes l)
			\iota(h_{(2)})
			(\beta\otimes l)^{-1} \\ &
			= \tilde\gamma(h_{(5)})\tilde\gamma^{-1}(h_{(1)}) \beta(h_{(2)})\beta^{-1}(h_{(4)})
			(\epsilon\otimes lh_{(3)}l^{-1})\\ &
			= 
			(\beta^{-1}\tilde\gamma)(h_{(3)})
			(\beta^{-1}\tilde\gamma)^{-1}(h_{(1)}) 
			\iota(lh_{(2)}l^{-1}).
		\end{align*}
		The injectivity of $\iota$ implies that $(\beta^{-1}\tilde \gamma, l)$ is a pair in involution for $H$.
	\end{proof}

	A Hopf algebra is called \emph{unimodular} if its distinguished group-like is equal to  its counit.
	
	\begin{lemma}\label{lem: UnimodPii}
		Suppose $H$ to be a generalised Taft algebra such that either $H$ or $H^\circ$ is unimodular. Then $H$ has a pair in involution.
	\end{lemma}
	\begin{proof}
		Assume without loss of generality $H^\circ$ to be unimodular.
		Let $N = 2^nj\geq 2$, with $j$ odd, $q$ a primitive $N$-th root of unity and $(a_{ij})$ a parameter matrix such that $H = H_q(a_{ij})$.
		By Lemma \ref{lem: integrals} the unimodularity of $H^\circ$ implies $a_{11} + a_{21} = 0$ modulo $N$.
		With respect to the matrices of powers and coefficients associated to $(a_{ij})$ this equation reads as
		$\mu_{11}2^{\mathfrak a_{11}} +  \mu_{21}2^{\mathfrak a_{21}} = 0$ modulo $2^n$. Thus, either 
		$\mu_{11} = \mu_{21} = 0$ or $\mathfrak a_{11} = \mathfrak a_{21}$. Applying Theorem \ref{thm: Classification} shows that $H$ has a pair in involution.		
	\end{proof}

	\begin{corollary}
		If the generalised Taft algebra $H$ is quasitriangular, it has a pair in involution.
	\end{corollary}

	\begin{proof}
		Let $N \geq 2$, $q$ a primitive $N$-th root of unity and $(a_{ij})$ a parameter matrix such that $H = H_q(a_{ij})$.
		Lemma \ref{lem: quasitriangular} asserts that $N$ is even and $a_{11} = a_{21} = \tfrac{N}{2}$. In particular, $a_{11} +a_{21} = 0$ modulo $N$.
		Lemma \ref{lem: integrals} shows that $H^\circ$ is unimodular and therefore $H$ has a pair in involution.
	\end{proof}

	\begin{remark}
		Drinfeld doubles are quasitriangular and unimodular, see Theorems 10.3.6 and 10.3.12 of \cite{Montgomery1993}, respectively.
		Applying Theorem \ref{thm: equivPiiHDualDouble} to the Hopf algebras of Lemma \ref{lem: Examples} shows that neither unimodularity nor quasitriangularity imply the existence of pairs in involution.	
	\end{remark}

	\section{Anti-Drinfeld doubles and pivotality}\label{sec: ADDandPiv}
	We investigate pairs in involution from a Morita theoretic viewpoint.
	This leads us to construct two generalised Taft algebras whose underlying algebras are isomorphic. Yet, only one of them admits a pair in involution, implying that such pairs are not a Morita equivalent property.
	We conclude this article by explaining the connection between pivotal categories and pairs in involution and commenting on a possible categorical description of such pairs. 
	
	\begin{lemma}\label{lem: IsoAlgNoPair}
		Let $q\in \k$ be a primitive $48$-th root of unity. The generalised Taft algebras  		
		$H \defeq H_q
		\left(\begin{smallmatrix}
			34 & 26\\
			4 & 4
		\end{smallmatrix}\right)$ and $L \defeq H_q\left(\begin{smallmatrix}
			34 & 26\\
			28 & 4
		\end{smallmatrix}\right)$ are isomorphic as algebras but only $L$ admits a pair in involution.
	\end{lemma}
	
	\begin{proof}
		One immediately verifies that $H$ and $L$ are generalised Taft algebras.
		Let $g,x,y \in H$ be the generators of $H$ and $\hat g , \hat x, \hat y$ the generators of $L$. The defining relations of the algebras $H$ and $L$ are
		\begin{gather*}
			x^{ \ord(34 \cdot 26 )} = x^{\ord( 20)} = 0, \qquad
			y^{\ord( 4 \cdot 4 )} = y^{\ord( 16 )} = 0, \qquad
			xy = q^{34 \cdot 4} yx = q^{40} yx,
			\\
			g^N = 1, \qquad 
			gx = q^{26} xg , \qquad
			gy = q^{4} yg,	
			\\
			\\
			\hat x^{\ord(34 \cdot 26 )} = \hat x^{\ord(20) } = 0, \qquad
			\hat y^{\ord(28 \cdot 4)} = \hat y^{\ord(16)} = 0, \qquad,
			\hat x \hat y = q^{34 \cdot 4} \hat y \hat x = q^{40} \hat y \hat x,
			\\
			\hat g^N = 1, \qquad 
			\hat g \hat x = q^{26} \hat x \hat g , \qquad
			\hat g \hat y = q^{4} \hat y \hat g.
		\end{gather*}
		Therefore, $H$ and $L$ are isomorphic as algebras.
		
		Note that $48 = 3 \cdot 16$. The parameter matrices of $H$ and $L$ modulo $16$ are		
		\begin{gather*}
			a^H_{ij} = 
			\begin{pmatrix}
				34 & 26 \\
				4 & 4	
			\end{pmatrix} =
			\begin{pmatrix}
				2 & 5 \cdot 2 \\
				4 & 4	
			\end{pmatrix}
			\text{ and }
			a^L_{ij} = 
			\begin{pmatrix}
				34 & 26 \\
				28 & 4	
			\end{pmatrix} =
			\begin{pmatrix}
				2 & 5 \cdot 2 \\
				3 \cdot 4 & 4	
			\end{pmatrix}.
		\end{gather*}
		In particular, we have 
		$( \mathfrak a_{ij}^H) = 
		(\mathfrak a_{ij}^L) = \left(\begin{smallmatrix}
			1 & 1\\
			2 & 2
		\end{smallmatrix}\right)$.
		The determinants and  powers of the matrices of coefficients are $\det (\mu_{ij}^H) = 12$ and $\det( \mu_{ij}^L) = 2$ as well as $\tau^H = 2$ and $\tau^L = 1$, respectively.
		Theorem \ref{thm: Classification} implies that $L$ has a pair in involution whereas $H$ does not.
	\end{proof}
	
	The next corollary follows readily from the fact that the class of generalised Taft algebras is closed under duality.
	
	\begin{corollary}
		There exist generalised Taft algebras $H$ and $L$ which are isomorphic as coalgebras such that only $L$ has a pair in involution.
	\end{corollary}
	
	\begin{remark}
		The generalised Taft algebras $H$ and $L$ of Lemma \ref{lem: IsoAlgNoPair} provide us with examples of Hopf algebras whose anti-Drinfeld doubles $A(H)$ and $A(L)$ are not Morita equivalent in the sense of Lemma \ref{lem: MoritaDfAntiDf}. 
		The Drinfeld and anti-Drinfeld double of $L$ are isomorphic as algebras by Theorem \ref{thm: PiiEquivIso}.
		The same argument as in the proof of Lemma \ref{lem: MoritaDfAntiDf}
		implies that $H$ would admit a pair in involution if $A(H)$ and $A(L) \cong D(L)$ were Morita equivalent in the above sense.
		This is a contradiction. 
	\end{remark}
	
	We end this article by characterising pairs in involution categorically.
	Recall that a \emph{pivotal category} is defined as a rigid monoidal category such that a natural monoidal isomorphism between every object and its bidual exists, see \cite[Definition 4.7.7]{Etingof2015}.
	This isomorphism can be encoded on the level of Hopf algebras, which is for example discussed in \cite{Andruskiewitsch2014}.
	\begin{definition}
		A \emph{pivotal Hopf algebra} is a pair $(H, \rho)$ of a Hopf algebra $H$ together with a group-like
		$\rho\in \Gr(H)$, its \emph{pivot}, such that
		$S^2(h) = \rho h \rho^{-1}$ for all $h \in H$.
	\end{definition}

	The following statement is a consequence of Equation \eqref{eq: PiiDDimpliesPiv}.
	\begin{lemma}
		A finite-dimensional Hopf algebra $H$ has a pair in involution if and only if $D(H)$ is pivotal.
	\end{lemma}
		
	The relationship between pivotal Hopf algebras and pivotal categories is well-known. For example, the first half of the following result is contained in the proof of \cite[Proposition 3.6]{Barrett1999}.
	
	\begin{lemma}
		The category $\LMod{H}_{\fin}$ of finite-dimensional left modules over a finite-dimensional Hopf algebra is pivotal if and only if $H$ is pivotal.
	\end{lemma}
	\begin{proof}
		Suppose $H$ to be pivotal with pivot $\rho \in \Gr(H)$.
		Given a finite-dimensional module $(M, \lact)$
		 over $H$ 
		we define
		$
			\omega_M \from M \rightarrow  \lbiDual{M},
			m \mapsto \can(\rho \lact m)
		$ 
		where $\lbiDual{M}$ denotes the left bidual of $M$ and `$\can$'  the canonical isomorphism of the underlying vector spaces.
		The definition of the dual module, see 
		\cite[III.6, Equation (6.5)]{Kassel1995}, shows that 
		$h \lact \can(m) = \can(S^2(h) \lact m)$ for every $h\in H$ and $m\in M$. Moreover,
		\begin{align*}
			\omega_{M} (h \lact m )& =
			\can(\rho h \lact m) =
			S^{-2}(\rho h) \lact \can (m)
			= \rho \rho^{-1} h \rho \lact \can (m)
			= h \lact \omega_{M}(m),
		\end{align*}
		for $h\in H$ and $m \in M$.
		Thus, $\omega_{M}$ is a morphism of modules whose inverse is given by $\omega_{M}^{-1} : \lbiDual{M} \rightarrow M, m \mapsto \rho^{-1} \lact \can^{-1}(m)$.
		The naturality of $\omega$ is verified by a straightforward computation and since $\rho$ is group-like $\omega$ is monoidal.

		Conversely, let $\LMod{H}_{\fin}$ be pivotal and   $\omega \from \Id \rightarrow  \lbiDual{(-)}$ a natural monoidal isomorphism.
		Consider $H$ as a module over itself with multiplication from the left as its action and define $\rho \defeq \can^{-1}(\omega_H(1)) \in H$.
		Since $\omega$ is monoidal $\rho$ is group-like and for all $h\in H$
		\begin{align*}
			S^2(h) \rho &= 
			S^2(h) \can^{-1}(\omega_H(1)) =
			\can^{-1}(h \lact \omega_H(1))
			\\
			&= \can^{-1}(\omega_H(h)) = 
			\can^{-1}(\omega_H(1 \cdot h)) =
			\can^{-1}(\omega_H(1))h =
			\rho h.
		\end{align*}
		 This proves $(H,\rho)$ to be a pivotal Hopf algebra.
	\end{proof}

	The \YetterDrinfeld modules over a finite-dimensional Hopf algebra $H$ are equivalent as a monoidal category to the modules over its Drinfeld double.
	By the two preceding lemmas, the \YetterDrinfeld modules over $H$ are pivotal if and only if $H$ admits a pair in involution.
	In conclusion there exists a connection between pairs in involution, anti-\YetterDrinfeld modules and pivotality.
	We would be interested in finding a more general characterisation of this interplay using an abstract notion of anti-\YetterDrinfeld modules as considered for example in \cite{Kobyzev2019}.
	
	\bibliographystyle{alpha}
	\bibliography{generalised_taft_algebras}

\newcommand{\etalchar}[1]{$^{#1}$}
\begin{thebibliography}{AAGI{\etalchar{+}}14}

\bibitem[AA17]{Andruskiewitsch2017}
Nicol\'{a}s Andruskiewitsch and Iv\'{a}n Angiono.
\newblock On finite dimensional {N}ichols algebras of diagonal type.
\newblock {\em Bull. Math. Sci.}, 7(3):353--573, 2017.

\bibitem[AAGI{\etalchar{+}}14]{Andruskiewitsch2014}
N.~Andruskiewitsch, I.~Angiono, A.~Garc\'{\i}a~Iglesias, B.~Torrecillas, and
  C.~Vay.
\newblock From {H}opf algebras to tensor categories.
\newblock In {\em Conformal field theories and tensor categories}, Math. Lect.
  Peking Univ., pages 1--31. Springer, Heidelberg, 2014.

\bibitem[AS98]{Andruskiewitsch1998a}
N.~Andruskiewitsch and H.-J. Schneider.
\newblock Lifting of quantum linear spaces and pointed {H}opf algebras of order
  {$p^3$}.
\newblock {\em Journal of Algebra}, 209(2):658--691, 1998.

\bibitem[AS02]{Andruskiewitsch2002}
Nicol\'{a}s Andruskiewitsch and Hans-J\"{u}rgen Schneider.
\newblock Pointed {H}opf algebras.
\newblock In {\em New directions in {H}opf algebras}, volume~43 of {\em Math.
  Sci. Res. Inst. Publ.}, pages 1--68. Cambridge Univ. Press, Cambridge, 2002.

\bibitem[BW99]{Barrett1999}
John~W. Barrett and Bruce~W. Westbury.
\newblock Spherical categories.
\newblock {\em Advances in Mathematics}, 143(2):357--375, 1999.

\bibitem[CM99]{Connes1999}
Alain Connes and Henri Moscovici.
\newblock {C}yclic cohomology and {H}opf algebras.
\newblock {\em Lett. Math. Phys.}, 48(1):97--108, 1999.

\bibitem[EGNO15]{Etingof2015}
Pavel Etingof, Shlomo Gelaki, Dmitri Nikshych, and Victor Ostrik.
\newblock {\em Tensor categories}, volume 205 of {\em Mathematical Surveys and
  Monographs}.
\newblock American Mathematical Society, Providence, RI, 2015.

\bibitem[Hec07]{Heckenberger2007}
Istv\'an Heckenberger.
\newblock {E}xamples of finite-dimensional rank 2 {N}ichols algebras of
  diagonal type.
\newblock {\em Compos. Math.}, 143(1):165--190, 2007.

\bibitem[Hec09]{Heckenberger2009}
Istv\'an Heckenberger.
\newblock {C}lassification of arithmetic root systems.
\newblock {\em Adv. Math.}, 220(1):59--124, 2009.

\bibitem[HK19]{Halbig2019}
Sebastian Halbig and Ulrich Kr\"ahmer.
\newblock A {H}opf algebra without modular pair in involution.
\newblock In {\em Geometric methods in physics {XXXVII}}, Trends Math., pages
  140--142. Birkh\"{a}user/Springer, Cham, 2019.

\bibitem[HKRS04a]{Hajac2004a}
Piotr~M. Hajac, Masoud Khalkhali, Bahram Rangipour, and Yorck Sommerh\"{a}user.
\newblock Hopf-cyclic homology and cohomology with coefficients.
\newblock {\em Comptes Rendus Math\'{e}matique. Acad\'{e}mie des Sciences.
  Paris}, 338(9):667--672, 2004.

\bibitem[HKRS04b]{Hajac2004}
Piotr~M. Hajac, Masoud Khalkhali, Bahram Rangipour, and Yorck Sommerh\"{a}user.
\newblock {S}table anti-{Y}etter--{D}rinfeld modules.
\newblock {\em C. R. Math. Acad. Sci. Paris}, 338(8):587--590, 2004.

\bibitem[HS10]{Hajac2010}
Piotr~M. Hajac and Yorck Sommerh\"{a}user.
\newblock {A}nti-{Y}etter--{D}rinfel’d {M}odules in {C}yclic {C}ohomology.
\newblock {\em privately communicated}, 2010.

\bibitem[Kas95]{Kassel1995}
Christian Kassel.
\newblock {\em Quantum groups}, volume 155 of {\em Graduate Texts in
  Mathematics}.
\newblock Springer-Verlag, New York, 1995.

\bibitem[KR93]{Kauffman1993}
Louis~H. Kauffman and David~E. Radford.
\newblock A necessary and sufficient condition for a finite-dimensional
  {D}rinfel'd double to be a ribbon {H}opf algebra.
\newblock {\em J. Algebra}, 159(1):98--114, 1993.

\bibitem[KS19]{Kobyzev2019}
Ivan Kobyzev and Ilya Shapiro.
\newblock A categorical approach to cyclic cohomology of quasi-{H}opf algebras
  and {H}opf algebroids.
\newblock {\em Applied Categorical Structures}, 27(1):85--109, 2019.

\bibitem[Mon93]{Montgomery1993}
Susan Montgomery.
\newblock {\em Hopf algebras and their actions on rings}, volume~82 of {\em
  CBMS Regional Conference Series in Mathematics}.
\newblock Published for the Conference Board of the Mathematical Sciences,
  Washington, DC; by the American Mathematical Society, Providence, RI, 1993.

\bibitem[Nen04]{Nenciu2004}
Adriana Nenciu.
\newblock Quasitriangular pointed {H}opf algebras constructed by {O}re
  extensions.
\newblock {\em Algebr. Represent. Theory}, 7(2):159--172, 2004.

\bibitem[Rad76]{Radford1976}
David~E. Radford.
\newblock The order of the antipode of a finite dimensional {H}opf algebra is
  finite.
\newblock {\em American Journal of Mathematics}, 98(2):333--355, 1976.

\bibitem[Rad12]{Radford2012}
David~E. Radford.
\newblock {\em Hopf algebras}, volume~49 of {\em Series on Knots and
  Everything}.
\newblock World Scientific Publishing Co. Pte. Ltd., Hackensack, NJ, 2012.

\end{thebibliography}
\end{document}